

\baselineskip=14pt
\parskip=10pt

\font\eightrm=cmr8 

\magnification=\magstephalf

\def\1{{\overline{1}}}
\def\2{{\overline{2}}}
\parindent=0pt
\overfullrule=0in

\def\frac#1#2{{#1 \over #2}}
\centerline
{\bf Enumerative Geometric Genealogy (Or: The Sex Life of Points and Lines)}
\bigskip
\centerline
{\it Shalosh B. EKHAD and Doron ZEILBERGER}
\bigskip
{\it ``One must be able to say at all times, instead of points and lines-tables and chairs''}-David Hilbert
\bigskip

{\bf Sexual Reproduction in the Fifth Planet of Star \#130103 of Galaxy \#4132}

The inhabitants of that planet, like in our own planet, are divided into two
(mutually exclusive!) {\it genders}, that, surprisingly, are also called
{\it male} and {\it female}. But the {\it sexual reproduction} is quite different!
In order to {\it be fruitful and multiply}, creatures of opposite sex can {\bf not} mate,
and while it is perfectly legal to have sex with the opposite gender, it could
never lead to progeny, and in the prehistory of that planet was 
considered ``sinful''.
If one wants to have children, then one must mate
with a member of the {\it same} gender, and then
the {\it baby} is {\bf always} of the {\bf opposite} gender to that of the parents.

[{\eightrm The details of the biology of the birth are open to speculations. Some scientists think
that each parent is pregnant with half of the baby, and the two halves are united after birth, 
while other think that the one who is going to be pregnant is chosen at random, and yet other people  believe that
there is no pregnancy, and the offspring is produced  instantaneously after the sexual act}.]

Another biological fact regarding the creatures in
that far-away planet is that any two individuals 
(that must be of the same gender, of course) can only have {\it one} child together.
In particular there can't be any twins, or triplets etc.
So no boy can have a full brother (of course, he can never have sisters!),
and no girl can have a full sister (of course, she can never have brothers!).
On the other hand, since {\it monogamy} is not  required (in fact it is forbidden, one must mate
with {\it all} eligible creatures that would give birth to a new creature), there are lots of {\it half-siblings}.

However, mating between half-siblings {\bf never} yields new creatures.
It turns out that if two individuals who share a parent (of course they can't share both parents, see above),
mate, their baby is an {\it exact clone} of that shared parent.

In our own planet Earth (the third planet of the star Sun in the Milky Way galaxy), as is
well-known, it all started with {\bf two} individuals, 
of {\it opposite} sex, Adam and Eve,
who were allowed to have numerous children, of both sexes. In order to keep going, 
there must have occurred some {\it incestuous} relationships (or we won't be here!)
[{\eightrm That's why incest only started to be sinful much later, after Lot and his daughters.}]

While two individuals (of the same sex, of course) may produce only {\it one} child, it sometimes happen that
two pairs of different parents give birth to {\bf identical} children.  Furthermore there is the following
{\it Law of Cogeny}.

{\bf The Law of Cogeny.} If the child of $A$ and $B$ is identical to the child of $A$ and $C$ , then
they are both identical to the child of $B$ and $C$. In that case the individuals $A,B,C$ are called
{\it cogenical}.

The obvious question  now is: 
{\it
``what is the minimal number of parentless {\it Adams}  God had to create {\it ab initio},
in order to account for the fact that there are now many {\it different} inhabitants?''}

[if a new child is born that is identical to an already existing creature, he or she is
expelled to another planet, no  clones are allowed.]

Unlike our own planet (Earth.Sun.MilkyWay@Universe.org), in that far away planet, two individuals (of the same sex, of course) do not suffice.
Let's call them $Adam_1$ and $Adam_2$. They are only allowed to have one daughter,
and that daughter  does not have any possible mating-mates!

But even with {\it three} starting individuals (once again, let them be male), let's call them
$$
Adam_1 \quad , \quad Adam_2 \quad , \quad Adam_3 \quad ,
$$
that species will not survive. Sure, the next (necessarily) female generation, now consists of
three females
$$
Eve_{1,2}:=DaughterOf(Adam_1,Adam_2) \quad , \quad
$$
$$
Eve_{1,3}:=DaughterOf(Adam_1,Adam_3) \quad , \quad
$$
$$
Eve_{2,3}:=DaughterOf(Adam_2,Adam_3) \quad ,
$$

but by the pigeon-hole principle, any two of these Eves share a parent, and
hence their matings only gives birth to clones of the Adams, and the species is doomed to 
only consist of six individuals: the three Adams and the three Eves.

But if the God of Planet \#5 of Star \#130103 of Galaxy \#4132
were kind enough to start with
{\bf four} Adams, (let's call the set of Adams `the $0^{th}$ generation'),
$$
Adam_1 \quad , \quad Adam_2 \quad , \quad Adam_3 \quad , \quad Adam_4 \quad ,
$$
then the population does increase.
Indeed, the first generation has now {\bf six} females
$$
Eve_{1,2}:=DaughterOf(Adam_1,Adam_2) \quad , \quad
Eve_{1,3}:=DaughterOf(Adam_1,Adam_3) \quad , \quad
$$
$$
Eve_{1,4}:=DaughterOf(Adam_1,Adam_4) \quad , \quad
Eve_{2,3}:=DaughterOf(Adam_2,Adam_3) \quad , \quad
$$
$$
Eve_{2,4}:=DaughterOf(Adam_2,Adam_4) \quad , \quad
Eve_{3,4}:=DaughterOf(Adam_3,Adam_4) \quad .
$$
The {\bf second} generation, (of males, of course) {\bf only} has three members (because half-sisters can't give birth to new creatures),
let's call them Abel, Cain, and Seth
$$
Abel:=SonOf(Eve_{1,2},Eve_{3,4}) \quad , \quad
Cain:=SonOf(Eve_{1,3},Eve_{2,4}) \quad , \quad
Seth:=SonOf(Eve_{1,4},Eve_{2,3}) \quad .
$$
Since it is perfectly OK to mate with a first-cousin, we now have,
at the  {\bf third} generation of creatures (alias the second generation of females),
three females, let's call them Sara, Rivka, and Lea.
$$
Sara:=DaughterOf(Abel,Cain) \quad , \quad
Rivka:=DaughterOf(Abel,Seth) \quad , \quad
Lea:=DaughterOf(Cain,Seth) \quad .
$$
Alas, now it seems that the poor inhabitants of our far-away planet are doomed to extinction, since any of the available
pairs $\{Sara, Rivka\}$, $\{Sara, Lea\}$, $\{Rivka, Lea\}$ are (half)-sisters, so can't give birth to new creatures.
Fortunately, the merciful God of that planet is allowing inter-generational mating!

Right now we have nine females! The six original Eves, and their three granddaughters. Of course if you
mate with your grandmother, the baby would be a clone of the man who is your father and her son, so in order
to give birth to new creatures,
$Eve_{1,2}$ is not allowed to know her granddaughters Sara and Rivka,
but, should mate with Lea. Similarly $Eve_{1,3}$ should only mate with  Rivka, etc.
Hence the second male generation (and the fourth altogether), consists of {\bf six} men,
let's call them Reuven, Shimon, Levi, Yehuda,Dan, and Naphtali:
$$
Reuven:=SonOf(Eve_{1,2},Lea) \quad , \quad Shimon:=SonOf(Eve_{1,3},Rivka) \quad , \quad Levi:=SonOf(Eve_{1,4},Sara) \quad ,
$$
$$
Yehuda:=SonOf(Eve_{2,3},Sara) \quad , \quad
Dan:=SonOf(Eve_{2,4},Rivka) \quad , \quad Naphtali:=SonOf(Eve_{3,4},Lea) \quad .
$$

How about the {\bf fifth} generation of creatures (alias the third generation of females)?
We currently have $4+3+6=13$ men alive. It is easy to see that they are $24$ additional fruitful
matings, but the great {\bf surprise} is that some of these new-born babies are identical.
It turns out that there are only $16$ {\bf different}-looking new babies.

{\bf Amazing Fact}: Reuven, Shimon, and Yehuda are cogenical.

In other words, the baby girl born to Reuven and Shimon is identical to the baby girl
born to Reuven and Yehuda, and both are the same as the baby girl born to Shimon and Yehuda.

By permuting the four Adams, we also have

$\bullet$ Reuven, Levi, and Dan are cogenical.

$\bullet$ Shimon, Levi, and Naphtali are cogenical.

$\bullet$ Yehuda, Dan, and Naphtali are cogenical.

The remaining twelve baby girls of the third female generation have no clones. Here they are:
$$
DaughterOf(Adam_1,Yehuda) \quad, \quad
DaughterOf(Adam_1,Dan) \quad, \quad
DaughterOf(Adam_1,Naphtali) \quad, \quad
$$
$$
DaughterOf(Adam_2,Shimon) \quad, \quad
DaughterOf(Adam_2,Levi) \quad, \quad
DaughterOf(Adam_2,Naphtali) \quad, \quad
$$
$$
DaughterOf(Adam_3,Reuven) \quad, \quad
DaughterOf(Adam_3,Levi) \quad, \quad
DaughterOf(Adam_3,Dan) \quad, \quad
$$
$$
DaughterOf(Adam_4,Reuven) \quad, \quad
DaughterOf(Adam_4,Shimon) \quad, \quad
DaughterOf(Adam_4,Yehuda) \quad .
$$

\vfill\eject

{\bf What If they were Five Adams}

In that case we get clones already at the fourth generation. In other words, there exist third-generation cogenical triples.

{\bf Four More Amazing Facts}: If there are five Adams, the following three females are cogenical:
$$
Eve_{1,2} \quad ,
$$
$$
DaughterOf(SonOf(Eve_{1,3}, Eve_{4,5}), SonOf(Eve_{2,4}, Eve_{3,5}) ) \quad ,
$$
$$
DaughterOf(SonOf(Eve_{1,4}, Eve_{3,5}), SonOf(Eve_{2,3}, Eve_{4,5}) ) \quad .
$$
As are the following three females:
$$
DaughterOf(Adam_1, SonOf(Eve_{2,3},Eve_{4,5})) \quad ,
$$
$$
DaughterOf(Adam_2, SonOf(Eve_{1,4},Eve_{3,5})) \quad ,
$$
$$
DaughterOf(Adam_5, SonOf(Eve_{1,3},Eve_{2,4})) \quad .
$$
As are the following three females:
$$
DaughterOf(Adam_1, SonOf(Eve_{2,3},Eve_{4,5})) \quad ,
$$
$$
DaughterOf(SonOf(E_{1,2} , E_{3,4}), SonOf(E_{1,3} , E_{2,5})) \quad , 
$$
$$
DaughterOf(SonOf(E_{1,4} , E_{2,5}), SonOf(E_{1,5} , E_{3,4})) \quad .
$$
As are the following three females:
$$
DaughterOf(SonOf(E_{1,2} , E_{3,4}), SonOf(E_{1,3} , E_{2,4})) \quad , 
$$
$$
DaughterOf(SonOf(E_{1,2} , E_{3,5}), SonOf(E_{1,3} , E_{2,5})) \quad , 
$$
$$
DaughterOf(SonOf(E_{2,4} , E_{3,5}), SonOf(E_{2,5} , E_{3,4})).
$$
Of course, there are numerous other such cogenical triples of females, obtained by permuting the Adams.

{\bf What About Six Adams}

With six Adams, we will not get look-alike babies any sooner, except if the six Adams all
live on the shore of the big ocean (that happens to be oval-shaped). Then we have

{\bf Yet another Amazing Fact}: If there are six Adams who live on the shore of the big oval ocean, then
the following second-generation creatures (i.e. first-generation men) are cogenical
$$
SonOf(Eve_{1,2},Eve_{3,4}) \quad, \quad 
SonOf(Eve_{1,5},Eve_{3,6}) \quad, \quad 
SonOf(Eve_{2,6},Eve_{4,5}) \quad.
$$
Of course, there are numerous other such cogenical triples of males, obtained by permuting the Adams.

{\bf Counting the Generations}

So far, in spite of the possibility of giving birth to identical babies, if they were four Adams, the
sizes of the successive generations (staring at generation $0$, the four Adams) are 
$$
4,6,3,3,6,16 \quad .
$$
It turns out that at the sixth-generation (i.e. third male generation) there are $84$ distinct new men, 
at the seventh generation, there are $1716$  distinct new women, 
while at the eighth generation, there are $ 719628$  distinct new men.

So the enumerating sequence starts with
$$
4, 6, 3, 3, 6, 16, 84, 1716, 719628 \dots \quad ,
$$
but we have no clue how it continues, or whether there is a `formula', or at least a polynomial-time algorithm
to compute the size of the $n$-th generation.

The enumeration of the successive generations with five Adams starts with:
$$
5, 10, 15, 90, 3495, \dots \quad ,
$$
and once again, we have no clue how it continues, while with six (general) Adams, it starts like this:
$$
6, 15, 45, 855, 342000 , \dots  \quad .
$$

{\bf Exegesis}

$\bullet$ Male $\rightarrow$ point \quad .

$\bullet$ Female $\rightarrow$ line \quad .

$\bullet$ The child of two males $\rightarrow$ the line joining two points \quad .

$\bullet$ The child of two females $\rightarrow$ the point of intersection of two lines \quad .

$\bullet$ {\bf Amazing fact} $\rightarrow$ Poncelet's result that the trilinear polar of a point, w.r.t. a triangle
exists [Wi1] (see also [We1]). Note that this is the simplest collinearity result that only uses iterations of the primitives
``point of intersection'' and ``line joining'' applied to four points in the plane in general position.

$\bullet$ {\bf Four more amazing facts} $\rightarrow$ the simplest concurrency theorems regarding five points in the plane in general
position, they are probably known, but who cares?
In some sense, once known, they are utterly trivial.
Does anyone care who was the first cave-woman that discovered that $7+2=4+5$?
We discovered (and immediately proved) them {\it ab initio}, using the Maple package {\tt GeometryMiracles} described below.

$\bullet$ {\bf Yet another amazing Fact} $\rightarrow$ Pascal's theorem ([Wi2],[We2]), where ``lying on the shore of the
oval ocean'' meant lying on a conic section.

{\bf Are there other Amazing Facts?}

Of course, once you have look-alike babies for non-obvious reasons that come from the original Adams,
there are lots of look-alike ``coincidences'' that are trivially implied by them, by
replacing  Adams by later-generation creatures, but are there any {\it genuinely} new ones?
In other words, to use a {\it fancy} word, can one find all the {\it syzygies}?
We don't know.

{\bf Enumerative Geometric Genealogy}

Let us briefly explain what we did.
At the $0^{th}$ generation, start with $k$ points on the plane, in {\it general position},
i.e. $k$ generic points $(s_1,t_1), \dots, (s_k,t_k)$ with $2k$ degrees of freedom. Then
to start a new odd generation (of lines), apply the operation
$$
Le((s,t),(s',t')):=[\, \frac{t-t'}{st'-s't} \, , \, \frac{s-s'}{st'-s't} \, ] \quad 
\eqno(Line)
$$
to all pairs of distinct already-existing points, getting
lots of new lines, (you may get some duplicates, but that's OK). 
Here $[a,b]$ is shorthand for the line $ax+by+1=0$.

To start a new even generation (of points), you do the analogous thing:
$$
Pt([s,t],[s',t']):=(\, \frac{t-t'}{st'-s't} \, , \, \frac{s-s'}{st'-s't} \, ) \quad .
\eqno(Point)
$$
Note that these formulas are identical, (hence we have {\it duality}).

Every object (point or line) has a {\it family tree}, and as we saw above, it is possible
for different family trees to yield the same object. Whenever that happens, we have a {\it miracle},
but most mathematicians call it a {\it theorem}.

{\bf How we  found out that we got scooped by Josh Cooper and Mark Walters}

Like all  self-respecting enumerators, we wanted to make sure that the above sequence
$$
4, 6, 3, 3, 6, 16, 84, 1716, 719628, \dots \quad ,
$$
of {\it new} points, and lines, at every generation, is brand new, and that no one 
considered this problem before. As we all know, the best way to do that
is to go to {\tt https://oeis.org/}. To our great delight, it was not there!
So we were positive that, {\it whp}, no one has ever considered this iteration of ``line joining''
and ``point of intersection'', starting with four general points, and enumerating the successive generations.
On June 17, 2014, one of us (DZ) had lunch with the guru, Neil Sloane, himself, and he claimed
that he {\it did} see it before. When we went back to his office, he quickly
came up with {\tt https://oeis.org/A140468} that counts the
number of creatures (of the same gender) born up to that generation, in other words, the sequence consisting
of the odd-indexed entries are the partial sums of our odd part, and ditto for the even-indexed entries.
Ah Well! The OEIS entry lead us to the interesting reference [CW] (that proved that if you start with
four points, the population explosion is {\bf doubly exponential}), as well as to the earlier beautiful article, [IR],
by Dan Ismailescu and Rado  Radoi{$\check c$}i\'c, that proved that the points are everywhere dense in the plane.
The analogous enumerations for five Adams and six Adams may be new (for what it's worth, we could not go very far).

{\bf Prohibiting Inter-Generational Matings}

We saw above that if you prohibit inter-generational matings, then starting
with four Adams would lead to extinction. However, with {\bf five} Adams,
we are safe. The (hopefully new) enumerating sequence  for this scenario is:
$$
5, 10, 15, 75, 2080, \dots \quad ,
$$
or with the convention of {\tt A140468}, it is
$$
5, 10, 20, 85, 2100, \dots , \quad .
$$

{\bf The Maple package GeometryMiracles}

Everything here was found thanks to the Maple package  {\tt GeometryMiracles}, available directly from 

{\tt http://www.math.rutgers.edu/\~{}zeilberg/tokhniot/GeometryMiracles}  \quad , \quad or via the front:

{\tt http://www.math.rutgers.edu/\~{}zeilberg/mamarim/mamarimhtml/egg.html} \quad ,

where there are numerous sample input and output files.

{\bf Further Work}

The Maple package {\tt GeometryMiracles} only iterates the primitives {\it Pt}, and {\it Le}
(that are essentially the same  operation), and looks for `surprises'. But one can  iterate
other primitives, like {\it MidPt}, that once again is a homosexual operation, and
its female counterpart {\it MidLe} (given by the same formula!), and more generally
{\it kPt} for any convex combination of the endpoints. Then we also have the
{\it heterosexual} operations of {\it Foot}, the projection of a point on a line,
and {\it Perp} the perpendicular line from a point to a line. 
There is also the {\it m\'enage \`a trois} operation, ``the circle passing through three distinct points'',
and its dual, ``the circle tangent to three distinct lines'', leading to yet another gender, {\it circle}, etc. etc.

Now one can start with just three points in general position, and easily discover {\it ab initio}, 
all the familiar theorems of triangle geometry, and  many new ones!

{\bf Encore I: A Schwartz-Tabachnikov Maple Rerun} 

The present \'etude in experimental mathematics was inspired by the
beautiful article, [ST],  of Richard Evan Schwartz and Serge Tabachnikov. The Maple package
{\tt RichardSerge} available directly from {\tt http://www.math.rutgers.edu/\~{}zeilberg/tokhniot/RichardSerge},\quad
or via the above front, independently confirms their pleasant surprises.

See the output files  {\tt http://www.math.rutgers.edu/\~{}zeilberg/tokhniot/oRichardSerge1} and \hfill\break
{\tt http://www.math.rutgers.edu/\~{}zeilberg/tokhniot/oRichardSerge2} \quad .

{\bf Encore II: Beyond Morley's Trisector Theorem} 

Euclid already knew that the angle-{\it bisectors}  of a triangle are {\it concurrent},
in other words, the side-lengths of the triangle formed by their intersections are all
zero (i.e. it degenerates into a point, the so-called   {\it incenter}).
More than $2100$ years later, Frank Morley (see [Wi3]) proved that the
analogous triangle for {\it trisectors} is no longer degenerate, but is {\it equilateral}. This triviality
fascinated many people, including such luminaries as Don Newman, John Conway, and
Alain Connes, who all published proofs.

But it took the genius of the first-named author of the present article (SBE), to come up with the next-in-line
theorem, a relation between the side-lengths of the triangle formed by the
{\it quadsectors}. For a degree-$14$ polynomial relating the {\it squares} of
the side-lengths, see

{\tt http://www.math.rutgers.edu/\~{}zeilberg/tokhniot/oGeneralizedMorley1} \quad .

{\bf Conclusion}

Roger Howe famously said (see [Z]):

{\it ``Everybody knows that mathematics is about miracles, but mathematicians have a special name for them: theorems''} \quad ,

but Doron Zeilberger (also see [Z]) retorted:

{\it ``Theorems are not miracles, but incestuous relationships between overdetermined inbred mathematical objects''} \quad .

Maybe they are both right!

{\bf Acknowledgment:}
Many thanks to  Neil Sloane, Eric Weisstein, and Paul Yiu for very useful pointers to the literature.

{\bf References}

[CR] Joshua Cooper and Mark Walters, {\it Iterated point-line configurations grow doubly-exponentially},
Disc. Comp. Geom. {\bf 43}(2010), 554-562, \quad {\tt http://arxiv.org/abs/0807.1549} \quad .

[IR] D. Ismailescu and R. Radoi{$\check c$}i\'c, {\it A dense planar point set from iterated line intersections},
Computational Geometry {\bf 27}(2004), 257-267.

[ST] Richard Evan Schwartz and Serge Tabachnikov, {\it Elementary Surprises in Projective Geometry},
{\tt http://arxiv.org/abs/0910.1952}. [A short version appeared in Math. Intell. {\bf 32}(2010), 31-34.]

[We1] Eric Weisstein, {\it Trilinear polar}, Wolfram Mathworld, \hfill\break
{\tt http://mathworld.wolfram.com/TrilinearPolar.html} \quad .

[We2] Eric Weisstein, {\it Pascal's theorem}, Wolfram Mathworld, \hfill\break
{\tt http://mathworld.wolfram.com/PascalsTheorem.html} \quad .

[Wi1] Wikipedia, {\it Trilinear polarity}, \hfill\break
{\tt http://en.wikipedia.org/wiki/Trilinear\_polarity}.

[Wi2] Wikipedia, {\it Pascal's theorem}, \hfill\break
{\tt  http://en.wikipedia.org/wiki/PascalTheorem } \quad .

[Wi3] Wikipedia, {\it Morley's trisector theorem}, \hfill\break
{\tt  http://en.wikipedia.org/wiki/Morley\_trisector\_theorem} \quad       .

[Z] Doron Zeilberger, {\it Opinion 44: Two lessons I learned from Shalosh B. Ekhad XIV's 2050 webbook} (posted Nov. 1, 2001), 
{\tt http://www.math.rutgers.edu/\~{}zeilberg/Opinion44.html} \quad .

\bigskip
\hrule
\bigskip
Shalosh B. Ekhad, c/o D. Zeilberger, Department of Mathematics, Rutgers University (New Brunswick), Hill Center-Busch Campus, 110 Frelinghuysen
Rd., Piscataway, NJ 08854-8019, USA.
\smallskip
Doron Zeilberger, Department of Mathematics, Rutgers University (New Brunswick), Hill Center-Busch Campus, 110 Frelinghuysen
Rd., Piscataway, NJ 08854-8019, USA. \hfill\break
Email: {\tt zeilberg at math dot rutgers dot edu} ;
url: {\tt http://www.math.rutgers.edu/\~{}zeilberg}.
\bigskip
\hrule
\bigskip
EXCLUSIVELY PUBLISHED IN \quad {http://www.math.rutgers.edu/\~{}zeilberg/pj.html} and {\tt arxiv.org}.
\bigskip
\hrule
\bigskip
June 19, 2014

\end